\documentclass{amsart}

\newcommand{\bdis}{\begin{displaymath}}
\newcommand{\edis}{\end{displaymath}}
\newcommand{\be}{\begin{equation}}
\newcommand{\ee}{\end{equation}}

\newcommand{\mcal}{\mathcal}

\newtheorem{theorem}{Theorem}%[section]

\theoremstyle{definition}

\newtheorem{cor}[theorem]{Corollary}

\theoremstyle{remark}

\begin{document}

\title{On Gram's law in the theory of the Riemann zeta function\footnote{paper published in ACTA ARITHMETICA, XXXII (1977), 107-113.}}

\author{Jan Mozer}

\address{Department of Mathematical Analysis and Numerical Mathematics, Comenius University, Mlynska Dolina M105, 842 48 Bratislava, SLOVAKIA}

\email{jan.mozer@fmph.uniba.sk}

%\keywords{Riemann zeta function, Friedmann cosmology}

%\begin{abstract}
%The aim of this paper is to show further results following those published in [5], and to relate the Riemann zeta function to the
%relativistic cosmology.
%\end{abstract}

\maketitle

{\bf 1.}\ Let $\{ t_\nu\}$ denote the sequence of the roots of the equation (see \cite{5}, p. 261):
\be
\vartheta(t)=\pi\nu ,
\ee
where (see \cite{5}, p. 383):
\be
\vartheta(t)=\frac{t}{2}\ln\left(\frac{t}{2\pi}\right)-\frac{t}{2}-\frac{\pi}{8}+\mcal{O}\left(\frac{1}{t}\right) ,
\ee
and $\nu$ runs over the set of all positive integers. Let us remind furthermore that (\cite{5}, p. 94) that
\be
Z(t)=e^{i\vartheta(t)}\zeta\left(\frac{1}{2}+it\right).
\ee

The Gram's law states that the zeroes of the function $\zeta\left(\frac{1}{2}+it\right)$ are separated by the points of the sequence $\{ t_\nu\}$ and
vice-versa. \\

E.C. Titchmarsh obtained in his work \cite{6} the following result
\be
\sum_{\nu=M+1}^N Z(t_\nu)Z(t_{\nu+1})\sim -2(c+1)N ,
\ee
where $M$ is a constant positive integer and $c$ is the Euler's constant. \\

Let $G(T)$ denote the number of intervals $(t_\nu,t_{\nu+1})\subset (0,T)$ containing a zero of the function $\zeta\left(\frac{1}{2}+it\right)$. E.C.
Titchmarsh deduced from the equation (4) that the following lower bound for $G(T)$:
\be
G(T)>A\frac{T^{2/3}}{\ln(T)}
\ee
holds true. \\

In this work, we show how it is possible to improve the results (4) and (5) in the sense of localization. Namely, we show the the following formula
holds true:
\be
\sum_{T\leq t_\nu\leq T+\bar{H}}Z(t_\nu)Z(t_{\nu+1})=-\frac{c+1}{\pi}\sqrt{T}\psi(T)\ln^2(T)+\mcal{O}(\sqrt{T}\ln^2(T)),
\ee
where $\psi(T)$ is an arbitrary increasing and not bounded from above function, and
\be
\bar{H}=\sqrt{T}\psi(T)\ln(T) .
\ee
Before formulating consequences of the formula (6), let us introduce some useful notations:
\be
\mu=\mu(T,\bar{H})=\max_{t\in [T,T+\bar{H}]}|Z(t)|,
\ee
where (see \cite{3}),
\be
|Z(t)|=\left|\zeta\left(\frac{1}{2}+it\right)\right| .
\ee
We call an interval $[\bar{t}_\nu,\bar{t}_{\nu+1}]$ regular if
\be
Z(\bar{t}_\nu)Z(\bar{t}_{\nu+1})<0 .
\ee
And finally, a regular interval contains an odd-number of zeroes of the function $Z(t)$, and therefore (with respect to (3)) also an odd-number of the
zeroes of the function $\zeta\left(\frac{1}{2}+it\right)$. \\

Let $G(T,\bar{H})$ denote the number of regular intervals which are in the interval $[T,T+\bar{H}]$. It follows from (6) and (8) (see \cite{6}, p. 105)
that
\begin{eqnarray}
& &
-\frac{c+1}{\pi}\sqrt{T}\psi(T)\ln^2(T)+\mcal{O}(\sqrt{T}\ln^2(T))\geq -\sum_{(\bar{t}_\nu)}|Z(\bar{t}_\nu)Z(\bar{t}_{\nu+1})|> \nonumber \\
& &
-\mu^2G(T,\bar{H}) ,
\end{eqnarray}
i.e. the following holds true:

\begin{cor}
\be
G(T,\bar{H})>A\frac{\sqrt{T}\psi(T)\ln^2(T)}{\mu^2}.
\ee
\end{cor}

By using results concerning the $\zeta\left(\frac{1}{2}+it\right)$ function (see \cite{5}, p. 116):
\be
\left|\zeta\left(\frac{1}{2}+it\right)\right|<A(\alpha,\beta)t^\alpha \ln^\beta(t),
\ee
we have as a consequence of (12) that

\begin{cor}
If (13) holds true then
\be
G(T,\bar{H})>A(\alpha,\beta)T^{1/2-2\alpha}\psi(T)\ln^{2-2\beta}(T) .
\ee
\end{cor}
The results of the type (13) listed in \cite{5}, p. 116 are to be extended by the following: \\
Hanke, \cite{1},
\be
\alpha=\frac{6}{37},\quad \beta=1;
\ee
Kolesnik, \cite{2},
\be
\alpha=\frac{173}{1067}, \quad \beta=\frac{331}{200} .
\ee
I am grateful to prof. Shitzel to draw my attention to the last result. \\

From (14) and (13) we have

\begin{cor}
\be
G(T,\bar{H})>AT^{721/2134}\psi(T)\ln^{-1.31}(T) .
\ee
\end{cor}

As well-known, the Lindeloef hypothesis states that for every $\epsilon>0$ we have
\be
\left|\zeta\left(\frac{1}{2}+it\right)\right|<A t^\epsilon,\quad t\geq T_0(\epsilon) .
\ee
And we have, combining (12) and (18):
\begin{cor}
If the Lindeloef hypothesis holds true then
\be
G(T,\bar{H})>A(\epsilon)T^{1/2-2\epsilon}\psi(T)\ln^2(T) .
\ee
\end{cor}
It follows from the Riemann hypothesis that the following estimate takes place (see \cite{5}, p. 350)
\be
\left|\zeta\left(\frac{1}{2}+it\right)\right|<\exp\left\{A\frac{\ln(t)}{\ln[\ln(t)]}\right\}=t^{\frac{A}{\ln[\ln(t)]}} .
\ee
And, combining (12) and (20), we conclude:
\begin{cor}
If the Riemann hypothesis holds true then
\be
G(T,\bar{H})>AT^{\frac{1}{2}-\frac{A}{\ln[\ln(T)]}}\psi(T)\ln^2(T) .
\ee
\end{cor}

In what follows, we give the proof of the formula (6).

{\bf 2.}\ In this section, we will discuss asymptotic terms related with relation (4). First of all we introduce the notation by \cite {5}, p. 98,
\be
g(t_\nu)=\sum_{m\leq \sqrt{t_\nu/2\pi}}\frac{\cos[t_\nu\ln(m)]}{\sqrt{m}} .
\ee
Now, one can find out by careful reading of pages 101-105 of the book \cite{5} that E.C. Titchmarsh has obtained the following relation
\be
\sum_{\nu=M+1}^N g(t_\nu)g(t_{\nu+1})=\frac{1}{2}(c+1)N+\mcal{O}\left(\frac{N}{\sqrt{\ln(N)}}\right)+
\mcal{O}(\sqrt{t_N}\ln^2(t_N)) .
\ee
The term
\bdis
\mcal{O}\left(\frac{N}{\sqrt{\ln(N)}}\right)
\edis
comes from (see \cite{6}, p. 102)
\be
\sum_{\nu=M+1}^N\left\{ \mcal{O}\left(\frac{1}{\sqrt{\ln(t_\nu)}}\right)\right\}=\mcal{O}\left(\frac{N}{\sqrt{\ln(N)}}\right) .
\ee
From (see \cite{5}, p. 261)
\be
Z(t_\nu)=2(-1)^\nu g(t_\nu)+\mcal{O}(t^{-1/4}_\nu),
\ee
and
\be
\sum_{\nu=M+1}^N\left\{ \mcal{O}(t^{-1/12}_\nu\ln(t_\nu))\right\}=\mcal{O}(t^{11/12}_N\ln(t_N)),
\ee
we deduce (see \cite{6}, p. 105) that
\be
\sum_{\nu=M+1}^N Z(t_\nu)Z(t_{\nu+1})=-4\sum_{\nu=M+1}^N g(t_\nu)g(t_{\nu+1})+\mcal{O}(t^{11/12}_N\ln(t_N)) .
\ee
And finally, (27) and (23) give us the result
\begin{eqnarray}
& &
\sum_{\nu=M+1}^NZ(t_\nu)Z(t_{\nu+1})= \nonumber \\
& &
-2(c+1)N+\mcal{O}\left(\frac{N}{\sqrt{\ln(N)}}\right)+\mcal{O}(t^{11/12}_N\ln(t_N))+\mcal{O}(\sqrt{t_N}\ln^2(t_N)) .
\end{eqnarray}
The last formula is a more accurate version of the relation (4). \\

Analyzing the Titchmarsh method - \cite{6}, p. 101-105 - we find that the following two tricks:
\begin{itemize}
\item[(a)] fixing the beginning of the sum:
\bdis
\sum_{\nu=M+1}^N \left( \dots \right)
\edis
i.e. setting $M$ to a constant, and
\item[(b)] a way in which the order of terms in the sum is inverted (this is, in fact, connected with the objection (a)) (see \cite{6}, p. 101)
\end{itemize}
are not fundamental in the Titchmarsh method of estimation of relevant quantities.  \\

Having this in mind we try to improve the relation (28) in the sense of the localization.

{\bf 3.}\ First of all, let
\be
H=\mcal{O}\{\sqrt{T}[\ln(T)]^k\},
\ee
where $k$ is a positive integer. Because of
\be
\sqrt{\frac{T+H}{2\pi}}-\sqrt{\frac{T}{2\pi}}=\frac{\frac{H}{2\pi}}{\sqrt{\frac{T+H}{2\pi}}+\sqrt{\frac{T}{2\pi}}}< A\frac{H}{\sqrt{T}}<
A[\ln(T)]^k ,
\ee
we have the inequality:
\be
\sum_{\sqrt{T/2\pi}\leq m\leq \sqrt{(T+H)/2\pi}}\frac{1}{\sqrt{m}}<A\frac{[\ln(T)]^k}{\sqrt[4]{T}} .
\ee
Subsequently, see (22), for $t_\nu\in [T,T+H]$ we have by (31) that
\be
g(t_\nu)=\sum_{m<\sqrt{T/2\pi}}\frac{\cos[t_\nu\ln(m)]}{\sqrt{m}}+\mcal{O}\left\{\frac{[\ln(T)]^k}{\sqrt[4]{T}}\right\}=
\sum_{(m)}+\sum^\prime .
\ee
Subsequently,
\begin{eqnarray}
P&=&\sum_{T\leq t_\nu<T+H}g(t_\nu)g(t_{\nu+1}) = \nonumber \\
&=&\sum_{T\leq t_\nu\leq T+H}\sum_{m\leq \sqrt{t_\nu/2\pi}}\frac{\cos[t_\nu\ln(m)]}{\sqrt{m}}\sum_{n\leq \sqrt{t_{\nu+1}/2\pi}}
\frac{\cos[t_{\nu+1}\ln(n)]}{\sqrt{n}}= \nonumber \\
&=&\sum_{(t_\nu)}\sum_{(m)}\sum_{(n)}+\sum_{(t\nu)}\sum^\prime\sum_{(n)}+\sum_{(t_\nu)}\sum_{(m)}\sum^{\prime\prime}+
\sum_{(t_\nu)}\sum^\prime\sum^{\prime\prime} .
\end{eqnarray}
At the end of this section we remind that
\be
\sum_{T\leq t_\nu\leq T+H}1=\frac{H}{2\pi}\ln\left(\frac{H}{2\pi}\right)+\mcal{O}\left(\frac{H^2}{T}\right) .
\ee

{\bf 4.}\ By using (25), (29), (31), (32), (34) and by the following inequality (see \cite{5}, p. 109)
\be
|Z(t)|=\left|\zeta\left(\frac{1}{2}+it\right)\right|<At^{1/6}\ln(t) ,
\ee
we have
\be
\sum_{(t\nu)}\sum^\prime\sum_{(n)}=\mcal{O}\left\{ H\ln(T)\frac{[\ln(T)]^k}{\sqrt[4]{T}}T^{1/6}\ln(T)\right\}=
\mcal{O}\left\{ T^{5/12}[\ln(T)]^{k+2}\right\} ,
\ee
\be
\sum_{(t_\nu)}\sum^\prime\sum^{\prime\prime}=\mcal{O}\left\{[\ln(T)]^{k+1}\right\} .
\ee
Also the following can be found in \cite{6}, p. 101:
\begin{eqnarray}
Q&=&\sum_{(t_\nu)}\sum_{(m)}\sum_{(n)}=\sum_{T\leq t_\nu\leq T+H}\sum_{m<\sqrt{T/2\pi}}\frac{\cos[t_\nu\ln(m)]}{\sqrt{m}}
\sum_{n<\sqrt{T/2\pi}}\frac{\cos[t_{\nu+1}\ln(n)]}{\sqrt{n}}=\nonumber \\
&=&\sum_{(m)}\frac{1}{\sqrt{m}}\sum_{(n)}\frac{1}{\sqrt{n}}\sum_{(t_\nu)}\cos[t_\nu\ln(m)]\cos[t_{\nu+1}\ln(n)]=\nonumber \\
&=&\frac{1}{2}\sum_{(m)}\frac{1}{\sqrt{m}}\sum_{(n)}\frac{1}{\sqrt{n}}\sum_{(t_\nu)}\cos[t_\nu\ln(m)-t_{\nu+1}\ln(n)]+\nonumber \\
& & +\frac{1}{2}\sum_{(m)}\frac{1}{\sqrt{m}}\sum_{(n)}\frac{1}{\sqrt{n}}\sum_{(t_\nu)}\cos[t_\nu\ln(m)+t_{\nu+1}\ln(n)]=
\frac{1}{2}\hat{\sum}+\frac{1}{2}\tilde{\sum} .
\end{eqnarray}
Now, the Titchmarsh method (\cite{6}, p. 101-105) can be used to estimate the quantity (38) with the result
\be
Q=\frac{1}{2}(c+1)\sum_{T\leq t_\nu\leq T+H}1+\mcal{O}\left\{ \sum_{T\leq t_\nu\leq T+H}\frac{1}{\sqrt{t_\nu}}\right\}+
\mcal{O}(\sqrt{T}\ln^2(T)) .
\ee
By (29) and (34) we have
\be
\sum_{T\leq t_\nu\leq T+H}1=\frac{H}{2\pi}\ln\left(\frac{T}{2\pi}\right)+\mcal{O}\{ [\ln(T)]^{2k}\},
\ee
and with respect to the formulae (39) and (40) we can write
\be
Q=\frac{c+1}{4\pi}H\ln\left(\frac{T}{2\pi}\right)+\mcal{O}(H\sqrt{\ln(T)})+\mcal{O}(\sqrt{T}\ln^2(T)) .
\ee
Now we can substitute (36), (37) and (41) into (33) to obtain
\begin{eqnarray}
P=\frac{c+1}{4\pi}H\ln\left(\frac{T}{2\pi}\right)&+&\mcal{O}(H\sqrt{\ln(T)})+\mcal{O}(\sqrt{T}\ln^2(T))+\nonumber \\
& &
+\mcal{O}\{T^{5/12}[\ln(T)]^{k+2}\}+\mcal{O}\{ [\ln(T)]^{k+1}\} .
\end{eqnarray}
And finally, we can put
\be
\bar{H}=\sqrt{T}\psi(T)\ln(T) ,
\ee
(which corresponds with $k=2$ in (29)) in order to obtain
\be
\sum_{T\leq t_\nu\leq T+H}g(t_\nu)g(t_{\nu+1})=\frac{c+1}{4\pi}\sqrt{T}\psi(T)\ln^2(T)+\mcal{O}(\sqrt{T}\ln^2(T)) .
\ee

{\bf 5.}\ The inequality
\be
|g(t_\nu)|<AT^{1/6}\ln(T),\quad t_\nu\in [T,T+\bar{H}] .
\ee
follows from (25) and (35). Furthermore, the formulae (25), (40), (43) and (45) imply the following:
\begin{eqnarray}
& &
\sum_{T\leq t_\nu\leq T+\bar{H}}Z(t_\nu)Z(t_{\nu+1})=\nonumber \\
& &
=\sum_{(t_\nu)}\{ 2(-1)^\nu g(t_\nu)+\mcal{O}(t_\nu^{-1/4})\}\{ 2(-1)^{\nu+1}g(t_{\nu+1})+\mcal{O}(t_{\nu+1}^{-1/4})\}=\nonumber \\
& &
=-4\sum_{(t_\nu)}g(t_\nu)g(t_{\nu+1})+\mcal{O}(\bar{H}\ln(T)T^{1/6}\ln(T)T^{-1/4})+\mcal{O}(T^{-1/2})=\nonumber \\
& &
=-4\sum_{(t_\nu)}g(t_\nu)g(t_{\nu+1})+\mcal{O}(T^{5/2}\psi(T)\ln^3(T)) .
\end{eqnarray}
And now, considering (44), we obtain the final result
\be
\sum_{T\leq t_\nu\leq T+\bar{H}}Z(t_\nu)Z(t_{\nu+1})=-\frac{c+1}{\pi}\sqrt{T}\psi(T)\ln^2(T)+\mcal{O}(\sqrt{T}\ln^2(T)) ,
\ee
which is exactly (6).

\end{document}